\documentclass{siamltex}

\pdfoutput=1
\usepackage[T1]{fontenc}
\usepackage[latin1]{inputenc}
\usepackage{fancybox}

\usepackage{amssymb}
\usepackage{amsmath}
\usepackage{mathpazo}

\usepackage{xcolor}
\usepackage{color}
\usepackage{float}
\usepackage[final]{graphicx}
\usepackage{enumerate}
\usepackage[%
  colorlinks%
  ,bookmarks={false}%
  ,pdftitle={paper3}%
  ,pdfsubject={Image Inpainting Based On Coherence Transport With Adapted Distance Functions}%
  ,pdfkeywords={}%
  ,pdfauthor={Thomas März <maerzt@ma.tum.de>}%
  ,plainpages={false}
  ]{hyperref}

\input{macros}

\title{Image Inpainting Based On Coherence Transport With Adapted Distance Functions}
\author{Thomas M\"{a}rz\thanks{Zentrum Mathematik,
        Technische Universit\"{a}t M\"{u}nchen,
        Boltzmannstr. 3, 85747 Garching, Germany
        ({\tt maerzt@ma.tum.de}). Manuscript as of \today.}}

\begin{document}

\maketitle

\begin{abstract} 
We discuss an extension of our method \emph{Image Inpainting Based on Coherence Transport}.
For the latter method the pixels of the inpainting domain have to be serialized into an ordered list.
Up till now, to induce the serialization we have used the distance to boundary map.
But there are inpainting problems where the distance to boundary serialization causes unsatisfactory inpainting results.
In the present work we demonstrate cases where we can resolve the difficulties by employing other distance functions
which better suit the problem at hand.
\end{abstract}



\section{Introduction}\label{Sect:Intro}
Non-texture image inpainting, also termed image interpolation, is the task of determining
the values of a digital image for a destroyed, or consciously masked, subregion of the image domain.

The simple idea of the generic single pass method -- which forms the basis for our method \emph{Image Inpainting Based on Coherence Transport} published in \cite{FBTM07}
 -- is to fill the inpainting domain by traversing its pixels in an onion peeling fashion 
from the boundary inwards and thereby setting new image values as weighted means of given or already calculated ones.
%The basic ingredients of the generic method are the specific serialization of the pixels which are to be inpainted and the choice of a weight for calculating means.

Telea has been the first to use such an algorithm in \cite{Telea04}: the pixels  are serialized according to their euclidean distance to the boundary of the inpainting domain
and the weight is such that image values are propagated mainly along the gradient of the distance map. By his choices of the weight and the pixel serialization, the method of Telea is not adapted to image.

\begin{figure}[t]
	\begin{center}
		\begin{minipage}{0.32\textwidth}
			\includegraphics[width=\textwidth]{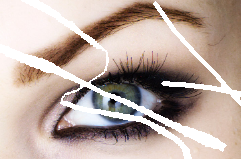}
		\end{minipage}
		\hfill
		\begin{minipage}{0.32\textwidth}
			\includegraphics[width=\textwidth]{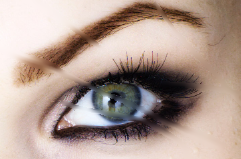}
		\end{minipage}
		\hfill
		\begin{minipage}{0.32\textwidth}
			\includegraphics[width=\textwidth]{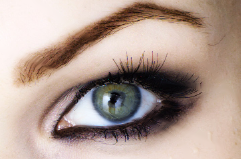}
		\end{minipage}
	\end{center}
	\caption{Scratch removal. Left: vandalized image (courtesy of Telea \cite[Fig 8.i]{Telea04}), the white scratches are the inpainting domain. 
	Middle: result of Telea's method. Right: result of our method}\label{Fig:EyeCompare}
\end{figure}

In \cite{FBTM07}, we addressed the adaption of the weight: our method uses an image dependent weight such that image values are propagated along the estimated tangents of
color lines which have been interrupted by the inpainting domain. That way we could improve the quality of the inpainting results compared to Telea (see figure \ref{Fig:EyeCompare}).
Beyond that, we have illustrated in \cite{FBTM07} that our method matches the high level of quality of the methods by \cite{Bert00}, \cite{Masnou98}, \cite{Masnou02}, and \cite{Tschump05}
while being considerably faster.

\begin{figure}[t]
	\begin{center}
		\begin{minipage}{0.32\textwidth}
			\includegraphics[width=\textwidth]{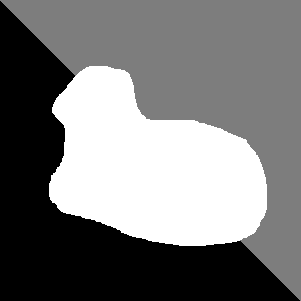}
		\end{minipage}
		\hfill
		\begin{minipage}{0.32\textwidth}
			\includegraphics[width=\textwidth]{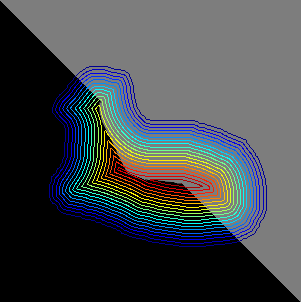}
		\end{minipage}
		\hfill
		\begin{minipage}{0.32\textwidth}
			\includegraphics[width=\textwidth]{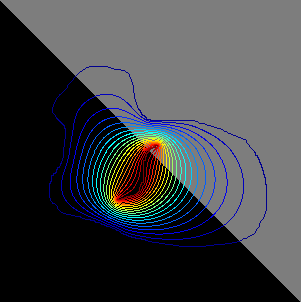}
		\end{minipage}
	\end{center}
	\caption{A broken diagonal ? Left: data image with inpainting domain in white. Middle: result of our method, pixel serialization by distance to boundary,
	levels of the distance to boundary map are overlaid. Right: result of our method, pixel serialization by another distance function.}\label{Fig:ABrokenDiag}
\end{figure}

But serializing the pixels by their distance to the boundary, as Telea and we in \cite{FBTM07} did, is not always a good idea as figure \ref{Fig:ABrokenDiag} illustrates.
The image in the middle shows the result which we obtain by using our method with distance to boundary serialization. 
The diagonal is not restored and the contours of the distance to boundary map indicate that this is due to the bad location 
of the skeleton which consists of the ridges of the distance map. Before reaching the skeleton, however, the performance is good and the diagonal is continued tangentially. 
The right image of figure \ref{Fig:ABrokenDiag} shows the result which we get if we use
another distance function (for pixel serialization) that suits the problem at hand better (see \S \ref{Sect:Harmonic}, figure \ref{Fig:BrokenDiagGoodSigmaCompare}).

The way of serializing the pixels is an important degree of freedom and in this paper we use this freedom to resolve some obstructions
which the distance to boundary serialization entails. 
The pixel serializations that we discuss are all induced by a more general type of distance function, so the \emph{new} parameter of the method
is in fact the user's choice of the distance function.
\medskip


\paragraph{Outline of the Paper}
In section \ref{Sect:Algo} we summarize the existing results and give a complete description of the basic algorithm.
In section \ref{Sect:DistFunc} we turn to concrete distance functions which are used to serialize the pixels of the inpainting domain.
We describe three approaches: distance by harmonic interpolation ({\S}\ref{Sect:Harmonic}), modified distance to boundary ({\S}\ref{Sect:ModDTB}),
and distance to skeleton ({\S}\ref{Sect:DistSkel}).
Section \ref{Sect:Conc} closes with a discussion.
\medskip

\section{Summary of Existing Results}\label{Sect:Algo}
%In this section we give a complete description of the algorithm and summarize the results of \cite{FBTM07}.
Our starting point is the generic algorithm for gray tone images. We assume that all gray tone images, seen as functions, take values in the real interval $\IV[0,{255}]$ and
we assume that gray tones are mapped onto $\IV[0,{255}]$ such that the natural order on the interval reflects  the order of the shades of gray by their brightness from
black ($=0$) to white ($=255$). Moreover, we will distinguish between discrete (digital) notions by using the index $h$ and continuous (analog) notions  where we omit the index.
Continuous notions are thought of as the high-resolution limit of the corresponding discrete ones. Finally, we identify pixels with their midpoints.
\medskip

\paragraph{Notation}
\begin{enumerate}[a)]
	\item $\Omega_{0,h}$ is the image domain, the matrix of pixels for the final, restored image $u_h: \Omega_{0,h} \to [0,1]$.
	\item $\Omega_h \subset \Omega_{0,h}$ is the inpainting domain whose values of $u_h$ have to be determined.
	\item $\Omega_{0,h} \wo \Omega_h$ is the data domain whose values of $u_h$ are given as $u_h|_{\Omega_{0,h} \wo \Omega_h}=u_{0,h}$.
	\item $\bd \Omega_h \subset \Omega_h$ is the discrete boundary, i.e., the set of inpainting pixels that have at least one neighbor in the data domain.
\end{enumerate}
Continuous quantities are defined correspondingly. Finally, we define discrete and continuous $\varepsilon$-neighborhoods by
\[
	B_{\varepsilon,h}(x) := \{y \in \Omega_{0,h} : |y-x| \leq \varepsilon\} \;, \quad B_{\varepsilon}(x) := \{y \in \Omega_{0} : |y-x| \leq \varepsilon\} \;.
\]
\medskip

\paragraph{Generic Algorithm}
The basic idea is to fill the inpainting domain from its boundary inwards by using weighted means of given or already calculated image values.
We assume  that the pixels which are to be inpainted have already been serialized, i.e. $\Omega_h = (x_1, x_2, \ldots, x_N)$ is an ordered list.
The following set
\[
	B^<_{\varepsilon,h}(x_k) := B_{\varepsilon,h}(x_k) \wo \{x_k, \ldots , x_N\} \;, \quad k=1,\ldots, N,
\]
denotes the neighborhood of the pixel $x_k$ consisting only of known or already inpainted pixels. Then, the algorithm reads as follows:
\begin{equation}\label{eqn:SinglePass}
	\begin{aligned}
		u_h|_{\Omega_{0,h} \wo \Omega_h} &=u_{0,h} \;,\\
		&\\
		u_h(x_k) &= \frac{\sum\limits_{y \in B^<_{\varepsilon,h}(x_k)} w(x_k,y) u_h(y)}{\sum\limits_{y \in B^<_{\varepsilon,h}(x_k)} w(x_k,y)} \;, \quad k=1,\ldots, N.
	\end{aligned}
\end{equation}
Here, $w(x,y) \geq 0$ are called the weights of the algorithm and we assume that
\[
	\sum\limits_{y \in B^<_{\varepsilon,h}(x)} w(x,y) >0 \;, \quad x \in \Omega_h \;.
\]
\medskip

\paragraph{Pixel Serialization}
In the generic algorithm, any serialization of the pixels can be used which geometrically goes from the boundary inwards in an onion peeling fashion. 
Now, let $T:\Omega \to \IV[0,{T_{\max}}]$ be a map with $T|_{\Gamma} = 0$, $\Gamma \subset \bd \Omega$, 
and which strictly grows into the interior of $\Omega$. We define the onion peels to be the level lines of $T$. 
Using such a map $T$ or rather a discretized version $T_h:\Omega_h \to \IV[0,{T_{\max}}]$, we serialize the pixels of $\Omega_h$ by 
\[
	T_h(x_j) < T_h(x_k) \quad \Rightarrow \quad j < k 
\]
into an ordered list $\Omega_h = (x_1, x_2, \ldots, x_N)$.

Clearly, there are different possible choices of a suitable $T$, but in the article \cite{FBTM07}
this degree of freedom has been fixed by setting $T=d$, where $d$ is the euclidean distance to the boundary $\bd \Omega$,   
\[
	d(x) = \dist(x,\bd \Omega) \;, \quad x \in \Omega \;.
\]
As pointed out in the introduction $T=d$ is not always a good choice and in the next section we will use generalized distances $T$ .
\medskip

\paragraph{Weights}
Here again, different choices for the weights are possible. We consider weights of the form
\begin{equation}\label{eqn:GenericW}
	w(x,y) = \frac{1}{|x-y|} \; k \left(x, (x-y) \cdot \varepsilon^{-1} \right) 
\end{equation}
with a kernel $k(x,\eta)$:	
\begin{enumerate}[a)]
	\item \emph{Normal Transport} (Telea's method). Telea, in \cite{Telea04}, uses the kernel 
		\[
			k(x,\eta) = \frac{|\SP< {N(x)}, {\eta} >|}{|\eta|} \;, \quad \mbox{with} \quad N(x) = \frac{\nabla T(x)}{|\nabla T(x)|}\;.
		\]
		Tough Telea has taken $T=d$, the kernel carries over to more general $T$.
	\item \emph{Coherence Transport} (our method of \cite{FBTM07}). For a given guidance field called $g$ we use the kernel
		\begin{equation}\label{eqn:CTkernel}
			k_{\mu}(x,\eta) = \sqrt{\frac{\pi}{2}} \mu \exp \left( -\frac{\mu^2}{2} \SP< {g^{\perp}(x)}, {\eta} >^2\right) \;.
		\end{equation}
		The guidance vector $g(x)$, supposed to be an approximate tangent, is computed by employing structure tensor analysis. 
		The set-up of the structure tensor $S_{\sigma,\rho}(x)$ is as follows:
		\begin{align*}
			\alpha_\sigma(y,x) &= \int\limits_{\Omega(x)} K_\sigma (y-h) \; dh \\
			v_\sigma (y,x) &= \frac{1}{\alpha_\sigma(y,x)} \cdot \int\limits_{\Omega(x)} K_\sigma (y-h) \cdot u(h) \; dh \\
			S_{\sigma,\rho}(x) &= \frac{1}{\alpha_\rho(x,x)} \cdot \int\limits_{\Omega(x)} K_\rho(x-y) \cdot \nabla_y v_{\sigma}(y,x) \cdot \nabla_y v^T_{\sigma}(y,x) \; dy
		\end{align*}
		where $\Omega(x) := \{y \in \Omega : T(y) \leq T(x)\} \cup \Omega_0 \wo \Omega$ and $K_\sigma$ is a Gaussian kernel. 
		The coherence vector $g(x)$ which enters $k_\mu$ as guidance is the eigenvector of $S_{\sigma,\rho}(x)$
		w.r.t. the minimal eigenvalue. Since $g$ depends on the image $u$, the weight does here, too.
		
		For more details on the usage of the structure tensor w.r.t. coherence transport inpainting see \cite{FBTM07} and \cite{MyDiss}.
		For other applications of structure tensor analysis and the related scale space theory see e.g. \cite{Aach06}, \cite{Weickert:DiffBook} and \cite{AGLM93}.
\end{enumerate}
\medskip

\paragraph{A Note on the Theory Behind} 
Because continuous theory is not an issue of this paper, we want at least make some comments.
For weights of the form \refEq{eqn:GenericW} we have shown in \cite{FBTM07} and \cite{MyDiss} that the high-resolution ($h \to 0$) vanishing-viscosity ($\varepsilon \to 0$) limit yields
transport/advection PDEs of the form:
\begin{enumerate}[a)]
	\item \emph{Normal Transport}. 
		\begin{align*}
			\SP< {N(x)} , {\nabla u(x)} >&= 0  \quad \mbox{in} \quad \Omega \wo \Sigma \;, \qquad & u|_{\bd \Omega} &= u_0|_{\bd \Omega} \;.
		\end{align*}
	\item \emph{Coherence Transport}. 
		\begin{align}\label{eqn:CTPDE}
			\SP< {c_\mu(x)} , {\nabla u(x)} >&= 0 \quad \mbox{in} \; \Omega \wo \Sigma \;, \; & u|_{\bd \Omega} &= u_0|_{\bd \Omega} \;,\;
			& \SP< {c_\mu(x)} , {N(x)} >& \geq \beta_\mu > 0\;.
		\end{align}
\end{enumerate}
In both cases the set $\Sigma$ denotes the set where $N(x) = \nabla T(x)/ |\nabla T(x)|$ is singular. The side condition $\SP< {c_\mu(x)} , {N(x)} > \geq \beta_\mu$
is the continuous counterpart of the pixel serialization: the function $T$ here plays the role of a Lyapunov function and this condition says that the characteristics
of the advection PDE evolve strictly into the interior of $\Omega$ and stop at $\Sigma$.
The general existence and well-posedness theory for such type of equations can be found in \cite{MyPaper1} and \cite{MyPaper2}.
This theory is applied to the analog inpainting model above in \cite[chapter 6]{MyDiss}.

Finally, in order to better understand the role of the parameter $\mu$ of the kernel $k_\mu$, we cite two results of \cite{FBTM07}.
Firstly, according to \cite[theorem 1]{FBTM07} every weight $w$ of the form \refEq{eqn:GenericW} corresponds to a transport field $c$.
Secondly, in \cite[theorem 2]{FBTM07}, when using the coherence transport kernel $k_\mu$, we have proved an asymptotic expansion of $c_{\mu}(x)$, w.r.t. $\mu \to \infty$, 
which implies the following limit behavior
\begin{equation}\label{eqn:MuLimit}
	\lim\limits_{\mu \to \infty} c_{\mu}(x) = 
	\begin{cases}
		g(x) & \;, \SP< {g(x)} , {N(x)}> > 0  \\
		-g(x) & \;, \SP< {g(x)} , {N(x)}> < 0  \\ 
		N(x) & \;, \SP< {g(x)} , {N(x)}> = 0  \\
	\end{cases} \quad .
\end{equation}
That means the coherence vector $g(x)$ guides in fact the coherence transport \refEq{eqn:CTPDE} and $\mu$ is the strength of this guidance.
\medskip

\paragraph{Interface of the Algorithm}
In the next section the algorithm will be performed only in its coherence transport version. 
That is, the weight function $w$ has the form of \refEq{eqn:GenericW} with the kernel $k_\mu$ given by equation \refEq{eqn:CTkernel}.
The execution of the coherence transport algorithm depends on the choice of four parameters:
\begin{itemize}
	\item $\varepsilon$, the averaging radius,
	\item $\mu$, the guidance strength of $k_\mu$,
	\item $\sigma$ and $\rho$, the scale parameters of the smoothing operations in the structure tensor $S_{\sigma,\rho}$.
\end{itemize}
Finally, the algorithm is supplied with the data image $u_0$ and a sorted list of the pixels which are to be inpainted. Any item of this list has the form
\[
	[ i \quad j \quad T_h(i,j)] \;,
\]
whereas $(i,j)$ are the pixel coordinates. % and $T_h(i,j)$ is the time value of the pixel. 
The list is sorted in ascending order of the $T_h(i,j)$ values.
\medskip

\section{Concrete Distance Functions}\label{Sect:DistFunc}
The generic algorithm of section \ref{Sect:Algo} depends on a prescribed pixel serialization, which orders the pixels from the boundary inwards.
In all previous experiments of \cite{FBTM07} the pixels were serialized by their euclidean distance to boundary.
The advantage in practice is that for all types of domains the approximate distance to boundary map is easy to generate 
by the Fast Marching Method (see \cite{Sethian99} and \cite{Kimmel04}).
But the disadvantage is that it is not always the best choice if one wants to get a nice looking inpainting result.
Here, we present three other ways of setting up a generalized discrete distance function $T_h$. 
We show a few synthetic examples where they yield better inpaintings than the distance to boundary.
\smallskip

\subsection{Distance by Harmonic Interpolation}\label{Sect:Harmonic}
\begin{figure}[tbp]
	\begin{center}
		\begin{minipage}{0.325\textwidth}
			\includegraphics[width=\textwidth]{figures/diag_tl2br_dist_mask.png}\\*[-1.5mm]
			{\scriptsize (a) damaged image} 
		\end{minipage}
		\hfill
		\begin{minipage}{0.325\textwidth}
			\includegraphics[width=\textwidth]{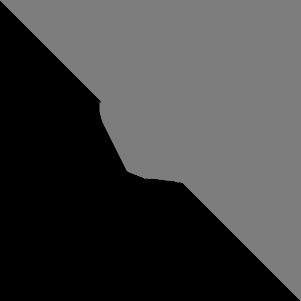}\\*[-1.5mm]
			{\scriptsize (b) inpainted image}
		\end{minipage}
		\hfill
		\begin{minipage}{0.325\textwidth}
			\includegraphics[width=\textwidth]{figures/diag_tl2br_dist_res_3_50_0p5_5_cont.png}\\*[-1.5mm]
			{\scriptsize (c) with contours of $d_h$}
		\end{minipage}
		\smallskip
				
		\begin{minipage}{0.325\textwidth}
			\includegraphics[width=\textwidth]{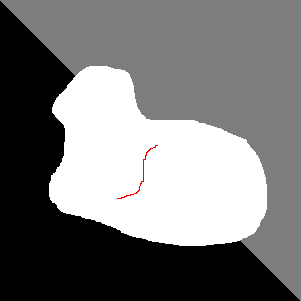}\\*[-1.5mm]
			{\scriptsize (d) damaged image} 
		\end{minipage}
		\hfill
		\begin{minipage}{0.325\textwidth}
			\includegraphics[width=\textwidth]{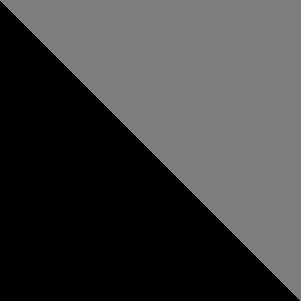}\\*[-1.5mm]
			{\scriptsize (e) inpainted image}
		\end{minipage}
		\hfill
		\begin{minipage}{0.325\textwidth}
			\includegraphics[width=\textwidth]{figures/diag_tl2br_goodsig_res_3_50_0p5_5_cont.png}\\*[-1.5mm]
			{\scriptsize (f) with contours of $T_h$}
		\end{minipage}
	\end{center}
	\caption{Broken diagonal. In the first row we have (a) the damaged image with the inpainting domain $\Omega_h$ in white,
	(b) the result inpainted using euclidean distance to boundary $d_h$, and (c) the result with contours of $d_h$.
	In the second row we have (d) the same damaged image with a specified stop set $\Gamma_h$ in red,
	(e) the result inpainted using distance by harmonic interpolation $T_h$, and (f) the result with contours of $T_h$.}\label{Fig:BrokenDiagGoodSigmaCompare}
\end{figure}

The broken diagonal is our first example. Figure \ref{Fig:BrokenDiagGoodSigmaCompare} (b) shows the result when using the distance to boundary $d_h$.
The desired inpainting result would be the restored diagonal but the diagonal is only partly continued correctly.
Figure \ref{Fig:BrokenDiagGoodSigmaCompare} (c) -- with the contours of $d_h$ overlaid -- demonstrates that the undesired effect is due to the \emph{wrong} location of the stop set. 
Figure \ref{Fig:BrokenDiagGoodSigmaCompare} (e) shows the result when using $T_h$ which we obtain by prescribing a better located stop set
and harmonic interpolation. Here, the algorithm is able to restore the diagonal.
For both cases we have used the same set of parameters $[\varepsilon, \mu, \sigma, \rho] = [3, 50, 0.5, 5]$.

Now, we describe the construction of the discrete distance function $T_h$.  Here, we can prescribe the location of the stop set arcs and their distance values.

Because the boundary is the start set, the discrete distance function $T_h$ shall equal zero on $\bd \Omega_h$.
In addition, we take at least one or more discrete curves $\Gamma_h^k$, $k \in \{1,\ldots,n\}$ , which are contained in $\Omega_h$
and should belong to the stop set $\Sigma_h$. Moreover, for every $\Gamma_h^k$ we specify a distance value $t_k > 0$.
The remainder of $T_h$ then, is calculated by harmonic interpolation. That is, we solve the discrete Laplace equation
\begin{align*}
	\Delta_h T_h &=0 \quad \mbox{in} \quad \Omega_h \wo \Sigma_h \;, & \Sigma_h &= \bigcup\limits_{k=1}^n \Gamma_h^k \;, \\
	T_h &= 0 \quad \mbox{on} \quad \bd \Omega_h \;, \qquad & T_h &= t_k \quad \mbox{on} \quad \Gamma_h^k \;, \quad k= 1,\ldots,n \;.
\end{align*}
Hereby, the discretization $\Delta_h$ of the Laplacian is due to the five-point-stencil
\begin{align*}
	\Delta_h T_h(i,j) &= T_h(i-1,j) + T_h(i,j-1) - 4 T_h(i,j) + T_h(i,j+1) + T_h(i+1,j) \;. 
\end{align*}

Since harmonic interpolation  provides a minimum and maximum principle, our construction of $T_h$ can be imagined as the
setting up a tent roof over  the domain $\Omega_h$ where  $\Gamma_h^k$  are the locations of the tent poles,
and every tent pole of $\Gamma_h^k$ has the length $t_k$. 

\begin{figure}[t]
	\begin{center}
	\begin{minipage}{0.32\textwidth}
		\includegraphics[width=\textwidth]{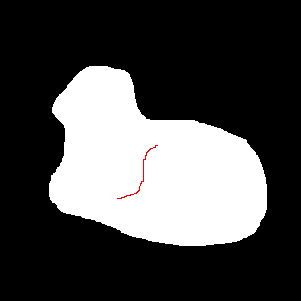}\\*[-1.5mm]
		{\scriptsize (a) white: $\Omega_h$, red: $\Gamma_h^1$ with $t_1 = 127$}
	\end{minipage}
	%\hfill
	\hspace{1cm}
	\begin{minipage}{0.32\textwidth}
		\includegraphics[width=\textwidth]{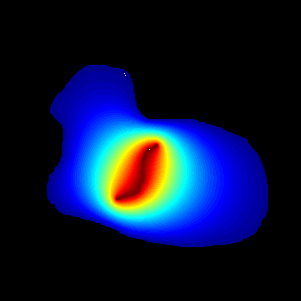}\\*[-1.5mm]
		{\scriptsize (b) contours of $T_h$ \phantom{$\Gamma_h^1$}}
	\end{minipage}
	\end{center}
	\caption{Single $\Gamma_h$ yields a valid $T_h$ }\label{Fig:VeryGoodTime}
\end{figure}

\begin{figure}[h]
	\begin{minipage}{0.32\textwidth}
		\includegraphics[width=\textwidth]{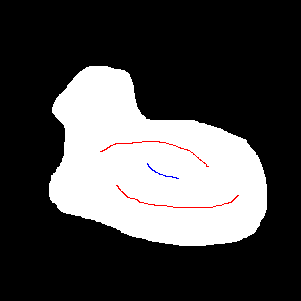}\\*[-1.5mm]
		{\scriptsize (a) white $\Omega_h$}
	\end{minipage}
	\hfill
	\begin{minipage}{0.32\textwidth}
		\includegraphics[width=\textwidth]{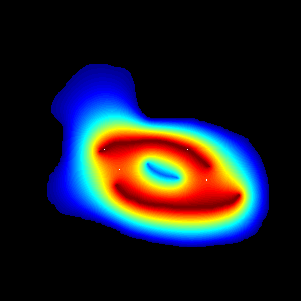}\\*[-1.5mm]
		{\scriptsize (b) non-valid $T_h$}
	\end{minipage}
	\hfill
	\begin{minipage}{0.32\textwidth}
		\includegraphics[width=\textwidth]{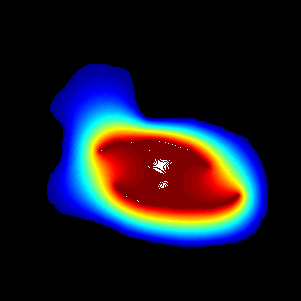}\\*[-1.5mm]
		{\scriptsize (c) valid $T_h$} 
	\end{minipage}
	\caption{(a) shows the domain $\Omega_h$ and the desired stop arcs, red: $\Gamma_h^1$ and $\Gamma_h^2$, blue: $\Gamma_h^3$.
        For (b) the distances on $\Gamma_h^k$ have been set to $t_1=t_2=250$, $t_3=50$ and we can see that the resulting $T_h$ has $\Gamma_h^3$ as local minimum and is thus not valid.
        For (c) the distances on $\Gamma_h^k$ have been set to $t_1=t_2=250$, $t_3=249$ and resulting $T_h$ is valid.}\label{Fig:BadTime}
\end{figure}

Unfortunately, not every choice of curves $\Gamma_h^k$, with time values $t_k$, results in a valid distance function: $T_h$ might have local minima
and so would not strictly increase into the interior of $\Omega_h$.
In the case of a single curve $\Gamma_h$ the resulting $T_h$ is a valid distance because of the minimum and maximum principle, see figure \ref{Fig:VeryGoodTime}. 
If there are two or more curves $\Gamma_h^k$, the question whether the resulting $T_h$ is valid or not
depends on the location of the curves in relation to each other and the differences of their prescribed values $t_k$  . 
Figure \ref{Fig:BadTime} (b) shows an example with three curves $\Gamma_h^1$, $\Gamma_h^2$ and $\Gamma_h^3$  with
\[
	t_1=t_2=250 \quad > \quad t_3=50 \;.
\] 
Here, all points of $\Gamma_h^3$ are local minima of $T_h$.
But, if we keep the geometry of $\Gamma_h^1$, $\Gamma_h^2$ and $\Gamma_h^3$, and change the prescribed distances to
\[
	t_1=t_2=250 \quad > \quad t_3=249 \;,
\] 
then $T_h$ does not have any minima (see figure \ref{Fig:BadTime} (c)). So, the resulting $T_h$ is admissible.
Generally speaking, if we have two or more curves $\Gamma_h^k$, with different prescribed distances $t_k$, and if the values $t_k$ are chosen unfavorably,
then the resulting $T_h$ might possess local minima on some of the $\Gamma_h^k$.

\begin{figure}[h]
	\begin{center}
		\begin{minipage}{0.325\textwidth}
			\includegraphics[width=\textwidth]{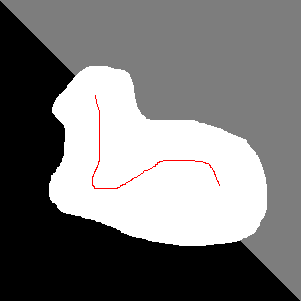}\\*[-1.5mm]
			{\scriptsize (a) damaged image}
		\end{minipage}
		\hfill
		\begin{minipage}{0.325\textwidth}
			\includegraphics[width=\textwidth]{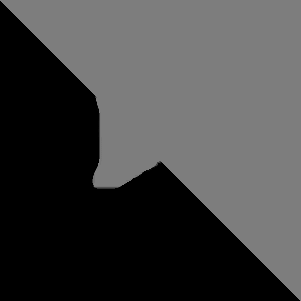}\\*[-1.5mm]
			{\scriptsize (b) inpainted image}
		\end{minipage}
		\hfill
		\begin{minipage}{0.325\textwidth}
			\includegraphics[width=\textwidth]{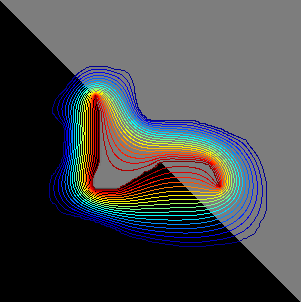}\\*[-1.5mm]
			{\scriptsize (c) with contours of $T_h$}
		\end{minipage}
	\end{center}
	\caption{Broken diagonal; method must fail if the stop set is deliberately placed badly. (a) the damaged image with the inpainting domain $\Omega_h$ in white with a single stop arc $\Gamma_h$ in red,
	(b) the result inpainted using distance by harmonic interpolation $T_h$, and (e) the result with contours of $T_h$. }\label{Fig:BrokenDiagBadSigma}
\end{figure}

Let us review our first example. Figure \ref{Fig:BrokenDiagGoodSigmaCompare} (e) illustrates that our $T_h$ -- which is that shown in figure \ref{Fig:VeryGoodTime} (b) -- works fine.
More generally, if we think of $\Sigma_h$ as an initial scratch, which has been dilated to all of $\Omega_h$ over the time $T_h$, then the backward
filling-in process, if $\Sigma_h$ is well located, makes the matching opposite sides come together.
Clearly, if we deliberately place $\Sigma_h$ badly, then the method must fail (see figure \ref{Fig:BrokenDiagBadSigma}, the parameters are $[\varepsilon, \mu, \sigma, \rho] = [3, 50, 0.5, 5]$, as before).

\begin{figure}[t]
	\begin{center}
		\begin{minipage}{0.325\textwidth}
			\includegraphics[width=\textwidth]{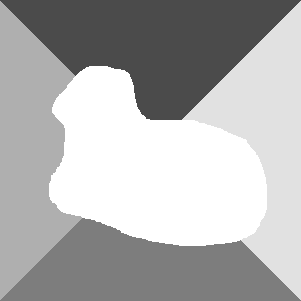}\\*[-1.5mm]
			{\scriptsize (a) damaged image}
		\end{minipage}
		\hfill
		\begin{minipage}{0.325\textwidth}
			\includegraphics[width=\textwidth]{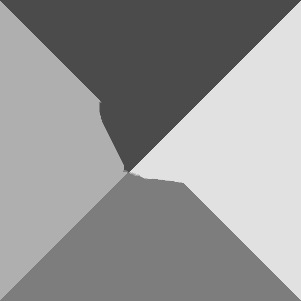}\\*[-1.5mm]
			{\scriptsize (b) inpainted image}
		\end{minipage}
		\hfill
		\begin{minipage}{0.325\textwidth}
			\includegraphics[width=\textwidth]{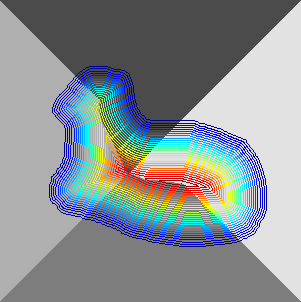}\\*[-1.5mm]
			{\scriptsize (c) with contours of $d_h$}
		\end{minipage}
		\smallskip
		
		\begin{minipage}{0.325\textwidth}
			\includegraphics[width=\textwidth]{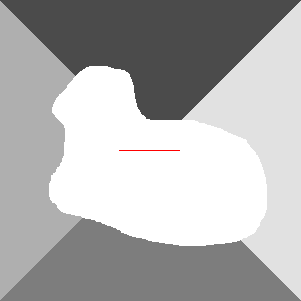}\\*[-1.5mm]
			{\scriptsize (d) damaged image}
		\end{minipage}
		\hfill
		\begin{minipage}{0.325\textwidth}
			\includegraphics[width=\textwidth]{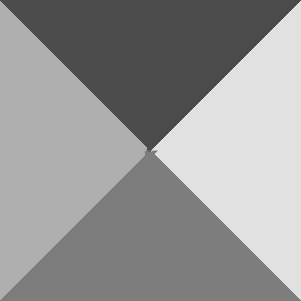}\\*[-1.5mm]
			{\scriptsize (e) inpainted image}
		\end{minipage}
		\hfill
		\begin{minipage}{0.325\textwidth}
			\includegraphics[width=\textwidth]{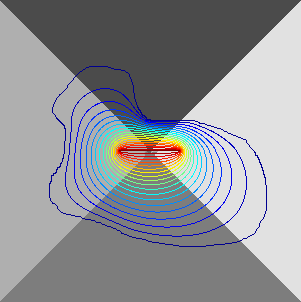}\\*[-1.5mm]
			{\scriptsize (f) with contours of $T_h$}
		\end{minipage}
	\end{center}
	\caption{Two broken diagonals. In the first row we have (a) the damaged image with the inpainting domain $\Omega_h$ in white,
	(b) the result inpainted using euclidean distance to boundary $d_h$, and (c) the result with contours of $d_h$.
	The result in (b) demonstrates that the undesired effect depends on how edges are located in relation to the inpainting domain.
	In the second row we have (d) the same damaged image with a specified stop set $\Gamma_h$ in red,
	(e) the result inpainted using distance by harmonic interpolation $T_h$, and (f) the result with contours of $T_h$.}\label{Fig:TwoBrokenDiagsGoodSigmaCompare}
\end{figure}

Now, we discuss two further examples.
Figure \ref{Fig:TwoBrokenDiagsGoodSigmaCompare} (a) shows two broken diagonals with the same inpainting domain as in figure \ref{Fig:BrokenDiagGoodSigmaCompare} (a).
We emphasize here that the appearance of an undesired effect depends on how the edge that needs to be continued is located in relation to the inpainting domain.
In figure \ref{Fig:TwoBrokenDiagsGoodSigmaCompare} (b) -- which is inpainted using distance to boundary -- the bottom-left-to-top-right diagonal is continued as desired, 
while the continuation of the top-left-to-bottom-right diagonal suffers from a badly located stop set.
In Figure \ref{Fig:TwoBrokenDiagsGoodSigmaCompare} (d) we have the damaged image with the single curve $\Gamma_h$, $t_1=127$ shown in red.
Figure \ref{Fig:TwoBrokenDiagsGoodSigmaCompare} (e) shows the result inpainted using $T_h$, the distance by harmonic interpolation.
Again, a good location of $\Sigma_h$ makes for a good result.
For both cases we have used the same set of parameters $[\varepsilon, \mu, \sigma, \rho] = [3, 50, 0.5, 5]$.

\begin{figure}[t]
	\begin{center}
		\begin{minipage}{0.45\textwidth}
			\includegraphics[width=\textwidth]{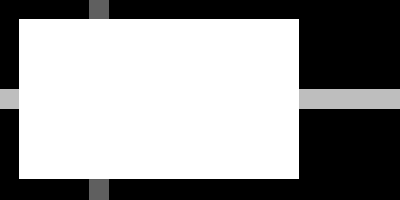}\\*[-1.5mm]
			{\scriptsize (a) damaged image}
		\end{minipage}
		\hfill
		\begin{minipage}{0.45\textwidth}
			\includegraphics[width=\textwidth]{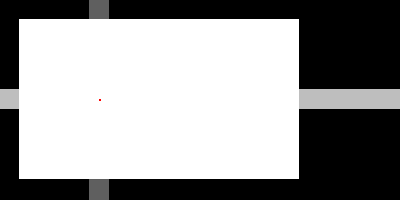}\\*[-1.5mm]
			{\scriptsize (d) damaged image, red $\Gamma_h$}
		\end{minipage}
		\smallskip
		
		\begin{minipage}{0.45\textwidth}
			\includegraphics[width=\textwidth]{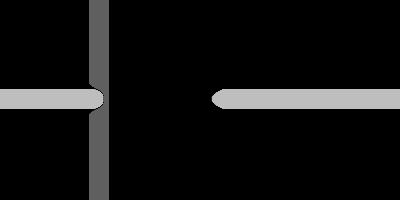}\\*[-1.5mm]
			{\scriptsize (b) inpainted image}
		\end{minipage}
		\hfill
		\begin{minipage}{0.45\textwidth}
			\includegraphics[width=\textwidth]{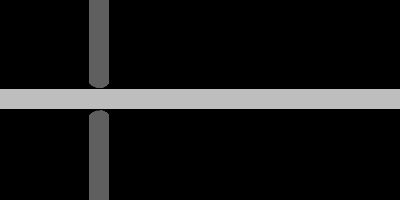}\\*[-1.5mm]
			{\scriptsize (e) inpainted image}
		\end{minipage}
		\smallskip
		
		\begin{minipage}{0.45\textwidth}
			\includegraphics[width=\textwidth]{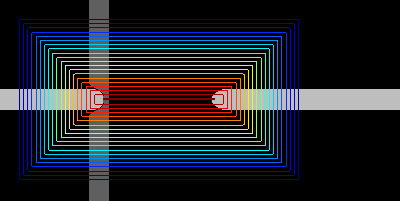}\\*[-1.5mm]
			{\scriptsize (c) with contours of $d_h$}
		\end{minipage}
		\hfill
		\begin{minipage}{0.45\textwidth}
			\includegraphics[width=\textwidth]{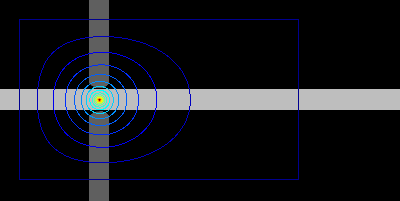}\\*[-1.5mm]
			{\scriptsize (f) with contours of $T_h$}
		\end{minipage}
	\end{center}
	\caption{Broken cross junction. In the first column we have (a) the damaged image with the inpainting domain $\Omega_h$ in white,
	(b) the result inpainted using euclidean distance to boundary $d_h$, and (c) the result with contours of $d_h$.
	In the second column we have (d) the same damaged image with the stop set $\Gamma_h$ being a single point in red (at the center of the cross junction),
	(e) the result inpainted using distance by harmonic interpolation $T_h$, and (f) the result with contours of $T_h$.}\label{Fig:BrokenJunctionGoodSigmaCompare}
\end{figure}

Figure \ref{Fig:BrokenJunctionGoodSigmaCompare} (a) shows a damaged cross junction.  A cross junction would, in any case, geometrically be
the simplest object for completion. Using the distance to boundary
we obtain the result shown in figure \ref{Fig:BrokenJunctionGoodSigmaCompare} (b).
Here, the stop set which is the central arc of the skeleton (see figure \ref{Fig:BrokenJunctionGoodSigmaCompare} (c)), has an unfavorable location because 
the bar coming from the right-hand side can never reach its counterpart.
Setting $\Gamma_h$ as the center of the cross junction (see figure \ref{Fig:BrokenJunctionGoodSigmaCompare} (d)) and using the distance by harmonic interpolation,
we are able to restore the cross junction (see figure \ref{Fig:BrokenJunctionGoodSigmaCompare} (e)).
Which of the bars is closed in the end depends on the coherence strength. The brighter bar has the higher
contrast w.r.t. the black background and is thus of stronger coherence. This is the reason why this bar is closed.
For both cases we have used the same set of parameters $[\varepsilon, \mu, \sigma, \rho] = [5, 100, 0.5, 10]$.

\subsection{Modified Distance to Boundary}\label{Sect:ModDTB}  
\begin{figure}[t]
	\begin{center}	
		\begin{minipage}{0.325\textwidth}
			\includegraphics[width=\textwidth]{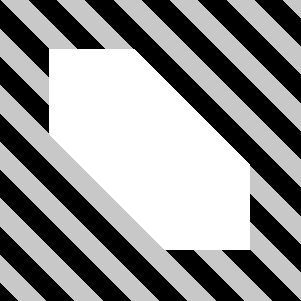}\\*[-1.5mm]
			{\scriptsize (a) damaged image}
		\end{minipage}
		\hfill
		\begin{minipage}{0.325\textwidth}
			\includegraphics[width=\textwidth]{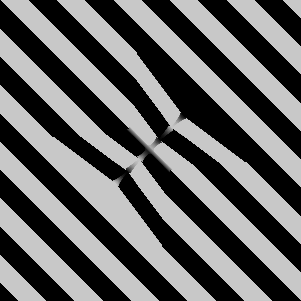}\\*[-1.5mm]
			{\scriptsize (b) inpainted image}
		\end{minipage}
		\hfill
		\begin{minipage}{0.325\textwidth}
			\includegraphics[width=\textwidth]{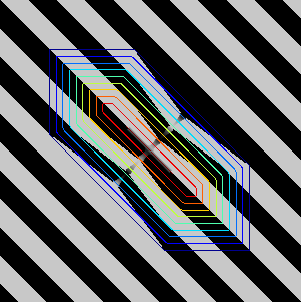}\\*[-1.5mm]
			{\scriptsize (c) with contours of $d_h$}
		\end{minipage}
		\smallskip
		
		\begin{minipage}{0.325\textwidth}
			\includegraphics[width=\textwidth]{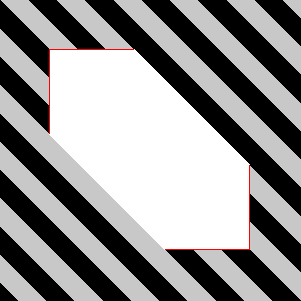}\\*[-1.5mm]
			{\scriptsize (d) red active boundary $\Gamma$}
		\end{minipage}
		\hfill
		\begin{minipage}{0.325\textwidth}
			\includegraphics[width=\textwidth]{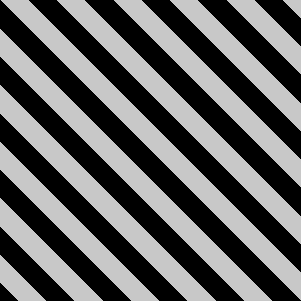}\\*[-1.5mm]
			{\scriptsize (e) inpainted image}
		\end{minipage}
		\hfill
		\begin{minipage}{0.325\textwidth}
			\includegraphics[width=\textwidth]{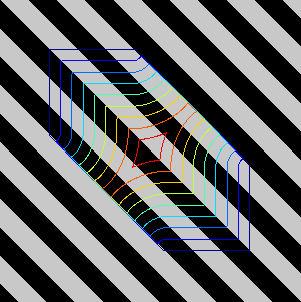}\\*[-1.5mm]
			{\scriptsize (f) with contours of $d_{*,h}$}
		\end{minipage}
	\end{center}
	\caption{Stripe pattern. The difficulty is that the damage is aligned with the pattern.
	In the first row we have (a) the damaged image with the inpainting domain $\Omega_h$ in white,
	(b) the result inpainted using euclidean distance to boundary $d_h$, and (c) the result with contours of $d_h$.
	In the second row we have (d) the same damaged image with the active part $\Gamma$ of the boundary in red,
	(e) the result inpainted using $d_{*,h}$ which is the euclidean distance to $\Gamma$, and (f) the result with contours of $d_{*,h}$.}\label{Fig:PlanarWavePBICompare}
\end{figure}

Figure \ref{Fig:PlanarWavePBICompare} (a) shows a damaged stripe pattern. Performing the algorithm with distance to boundary and the set of parameters
$[\varepsilon, \mu, \sigma, \rho] = [5, 100, 0.5, 10]$ yields the result shown in figure \ref{Fig:PlanarWavePBICompare} (b). 
The difficulty, is that the coherence vector which is tangent to the edges is orthogonal to the lower left and the upper right segment of $\bd \Omega_h$. 
Thus, the transport vector $c$ switches to the normal $N$ as in the exceptional case of equation \ref{eqn:MuLimit}.
Figure \ref{Fig:PlanarWavePBICompare} (c) confirms the switch of the transport to $N$.

To combat this we suggest a modification of the distance map set-up.
The euclidean distance to boundary map $d$ is the viscosity solution of
\begin{align*}
	|\nabla d| &= 1 \quad \mbox{in} \quad \Omega \;, & d|_{\bd \Omega} = 0 \;. 
\end{align*}
The modification is that we search for the euclidean distance $d_*$ to a subset $\Gamma$ of the boundary $\bd \Omega$, i.e.,
\begin{align*}
	|\nabla d_*| &= 1 \quad \mbox{in} \quad \Omega \;, & d_*|_{ \Gamma} = 0 \;. 
\end{align*}

We classify the points which shall not belong to $\Gamma$.
Assume $x \in \bd \Omega$ satisfies $d_*(x) = 0$, then the boundary normal is given by $N(x) = \nabla d_*(x)$.
Now, if we have at $x$
\[
	 \SP< {g(x)}, {N(x)}>^2 = 0 \;,
\]
whereas $g$ is the guidance vector, then either there is no guidance ($|g(x)| = 0$)
or the guidance vector does not point inwards. Such a boundary point shall not belong to $\Gamma$.
In fact, we use the stronger condition
\[
	0 \leq \SP< {g(x)}, {N(x)}>^2 \leq \gamma
\]
with a threshold parameter $0 < \gamma \leq 1$. Complementarily, the set of active boundary points $\Gamma$ is given by
\[
	\Gamma = \{x \in \bd \Omega: \SP< {g(x)}, {N(x)}>^2 > \gamma \} \;.
\]
Clearly, the new parameter $\gamma$ must be chosen such that $\Gamma$ is not empty.

We have applied this idea to the stripe pattern: the red lines in figure \ref{Fig:PlanarWavePBICompare} (d) are the active boundary points and the result
is shown in figure \ref{Fig:PlanarWavePBICompare} (e).
Our standard parameters are $[\varepsilon, \mu, \sigma, \rho]=[5, 100, 0.5, 10]$, while the additional parameter is set to $\gamma = 0.1$.
The discrete approximation $d_{*,h}$ was computed using the fast marching method.
The overlaid contour plot of $d_{*,h}$ in figure \ref{Fig:PlanarWavePBICompare} (f) shows that the guidance vector always points inwards.

\subsection{Distance to Skeleton}\label{Sect:DistSkel}  
The third approach to obtaining a serialization is to use the distance to a prescribed stop part of the skeleton. Let $\Sk^k$, $k \in \{1,\ldots,n\}$ be curves
in the image domain $\Omega_{0}$.  Those curves $\Sk^k$ which are contained in $\Omega$ will later belong to the skeleton $\Sk$. Let, then, $T_*$ be the viscosity solution of
\begin{align*}
	|\nabla T_{*}| &= 1 \quad \mbox{in} \quad \Omega_0 \;, &
	T_{*} &= 0 \quad \mbox{on} \quad \Sk^k \;, \quad k \in \{1,\ldots,n\} \;,
\end{align*}
and let
\[
	T_{*,\max} = \max\limits_{x \in \Omega}T_*(x) \;.
\]
The desired distance function is defined by
\[
	T(x) = T_{*,\max} - T_*(x) \;, \quad x \in \Omega \;.
\]
Warning: as in the case of harmonic interpolation (see \S \ref{Sect:Harmonic}) one must check if $T$ is admissible, i.e., if $T$ is free of local minima.

\begin{figure}[t]
	\begin{center}
		\begin{minipage}{0.325\textwidth}
			\includegraphics[width=\textwidth]{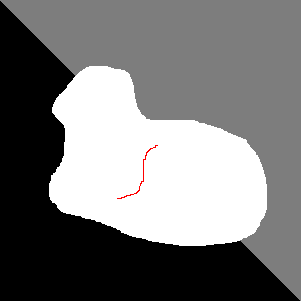}\\*[-1.5mm]
			{\scriptsize (a) damaged image}
		\end{minipage}
		\hfill
		\begin{minipage}{0.325\textwidth}
			\includegraphics[width=\textwidth]{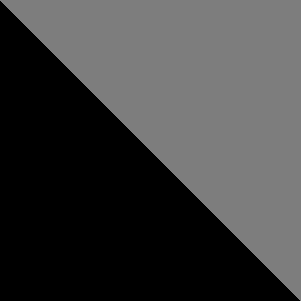}\\*[-1.5mm]
			{\scriptsize (b) inpainted image}
		\end{minipage}
		\hfill
		\begin{minipage}{0.325\textwidth}
			\includegraphics[width=\textwidth]{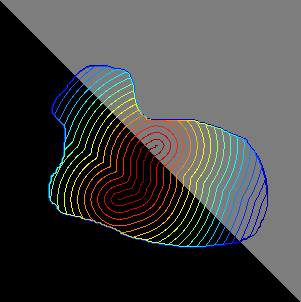}\\*[-1.5mm]
			{\scriptsize (c) with contours of $T_h$}
		\end{minipage}
	\end{center}
	\caption{Broken diagonal. (a) shows the damaged image, white $\Omega_h$, with the prescribed skeleton arc $\Sk_h^1$ in red.
	(b) is the result inpainted using distance to skeleton $T_h$. (c) is the result with contours of $T_h$.}\label{Fig:BrokenDiagDistSigma}
\end{figure}

Figure \ref{Fig:BrokenDiagDistSigma} shows the result for the example of the broken diagonal. 
The red curve in figure \ref{Fig:BrokenDiagDistSigma} (a) defines $\Sk_h^1$ (discrete) which shall belong to the skeleton.
The discrete approximation $T_{*,h}$ was computed using the fast marching method and the set of parameters is $[\varepsilon, \mu, \sigma, \rho] = [3, 50, 0.5, 5]$,
the same as we used for figure \ref{Fig:BrokenDiagGoodSigmaCompare} (e). 
The inpainted result here (figure \ref{Fig:BrokenDiagDistSigma} (b)) is the same as in figure \ref{Fig:BrokenDiagGoodSigmaCompare} (e), 
but the distance function has changed (compare figures \ref{Fig:BrokenDiagDistSigma} (c) and \ref{Fig:BrokenDiagGoodSigmaCompare} (f)). 
If there is only one curve $\Sk_h^1$ which is completely contained in $\Omega$, then distance by harmonic interpolation
(with $\Gamma_h^1 = \Sk_h^1$) and distance to skeleton will yield very similar results.
If there are two or more curves, then distance by harmonic interpolation allows for
different distance values on $\Gamma_h^k$, while distance to skeleton has exactly one distance value on all of the curves $\Sk_h^k$.
In contrast to distance to skeleton, the distance by harmonic interpolation method requires $T_h|_{\bd \Omega_h} = 0$.
 
\begin{figure}[t]
	\begin{center}
		\begin{minipage}{0.325\textwidth}
			\includegraphics[width=\textwidth]{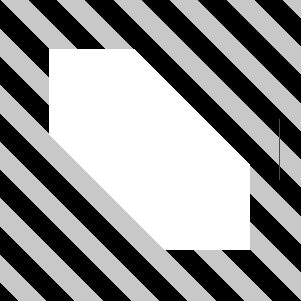}\\*[-1.5mm]
			{\scriptsize (a) damaged image}
		\end{minipage}
		\hfill
		\begin{minipage}{0.325\textwidth}
			\includegraphics[width=\textwidth]{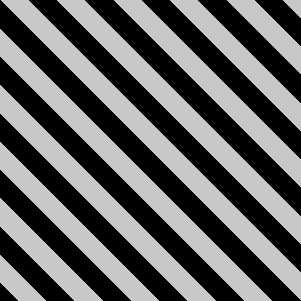}\\*[-1.5mm]
			{\scriptsize (b) inpainted image}
		\end{minipage}
		\hfill
		\begin{minipage}{0.325\textwidth}
			\includegraphics[width=\textwidth]{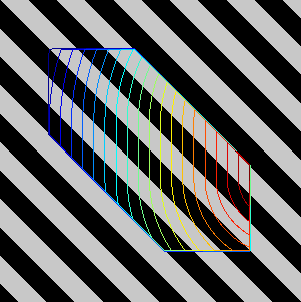}\\*[-1.5mm]
			{\scriptsize (c) with contours of $T_h$}
		\end{minipage}
	\end{center}
	\caption{Stripe pattern. (a) shows the damaged image, white $\Omega_h$, with the red arc $\Sk_h^1$ outside $\Omega_h$.
	(b) is the result inpainted using $T_h$, the distance to ''skeleton'' $\Sk_h^1$. (c) is the result with contours of $T_h$.}\label{Fig:PlanarWaveDistSigma}
\end{figure}

Moreover, since $T_*$ is defined on $\Omega_0$, we can place $\Sk^k$ outside of the inpainting domain $\Omega$.
We use this possibility to restore the stripe pattern in figure \ref{Fig:PlanarWaveDistSigma}.
The parameters are $[\varepsilon, \mu, \sigma, \rho]=[5, 100, 0.5, 10]$, the same as we used for figure \ref{Fig:PlanarWavePBICompare}.
It is obvious from the level lines of $T_h$ (see \ref{Fig:PlanarWaveDistSigma} (c)) that the guidance vector always points inwards.
Thus, we get the desired result again here.

Remark: The construction here, as well as that of \S \ref{Sect:ModDTB} produces distance functions whereas not the whole boundary belongs
to the start set. Besides that, it is possible that parts of the boundary belong in fact to the stop set. 
\medskip

\section{Discussion}\label{Sect:Conc}  
In the previous section we have demonstrated that there are inpainting problems where a distance to boundary induced serialization of the pixels
is not favorable and we have constructed other distances which are able to resolve the difficulties which appear. 

\begin{figure}[t]
	\begin{center}
		\begin{minipage}{0.475\textwidth}
			\includegraphics[width=\textwidth]{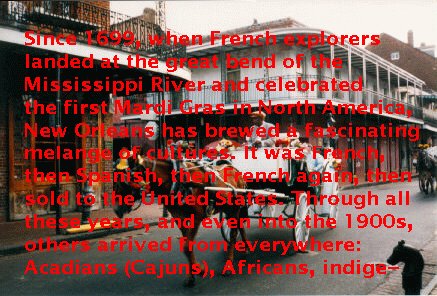}\\*[-1.5mm]
			{\scriptsize (a) original image (courtesy of \cite[figure 6]{Bert00})}
		\end{minipage}
		\hfill
		\begin{minipage}{0.475\textwidth}
			\includegraphics[width=\textwidth]{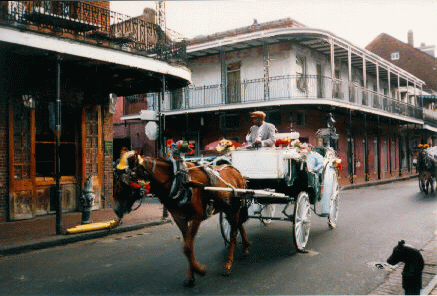}\\*[-1.5mm]
			{\scriptsize (b) inpainted: $[\varepsilon, \mu, \sigma, \rho] = [4, 25, 2, 3]$}
		\end{minipage}
	\end{center}
	\caption{ Removal of superimposed text. The letters in (a) are the inpainting domain $\Omega_h$}\label{Fig:NewOrleans}
\end{figure}

So far we have considered synthetic images, because their image geometry is easy to understand. 
Thus, we were able to construct distances which are adapted to the image or rather to an expected result.
When we face natural inpainting problems, as shown in figures \ref{Fig:EyeCompare} (left) and \ref{Fig:NewOrleans} (a), 
it is not as easy to set up an adapted distance function.
This is because
\begin{itemize}
	\item the geometry of the image is harder to understand,
	\item the damaged region is complicated.
\end{itemize}
Moreover, if the damaged region $\Omega$ consists of many connected components and one wants to use our distance by harmonic interpolation technique
(\S \ref{Sect:Harmonic} ), one must prescribe arcs $\Gamma_h^k$ for every single component of connectivity;
for the example of figure \ref{Fig:NewOrleans} (a) that means for every single letter.
This can be time consuming (for the user, not for the computer).

In contrast, the distance to boundary map can be computed fast and easily for every type of inpainting domain.
And, inpainting using distance to boundary often produces results of high quality 
when applied to natural inpainting problems (see figures \ref{Fig:EyeCompare} (right) and \ref{Fig:NewOrleans} (b)).
This is because for natural inpainting problems the damage often is such that
color lines have been broken by scratches (By scratches we mean rather thin and lengthy damages).
If the damage is of this type, the skeleton of $\Omega$, being a simplified version of the scratch, is well placed
and so the filling-in process makes the matching opposite sides come together.
The images shown in figures \ref{Fig:EyeCompare} (right) and \ref{Fig:NewOrleans} (a), like many other natural inpainting problems, have this type of damage.
Thus, our inpainting method using distance to boundary is able to produce results pleasing to the eye.

\section{Acknowledgements}
The author would like to thank Folkmar Bornemann for his advice and the inspiring discussions.

\bibliographystyle{siam}
\bibliography{DPaper3}

\end{document}